\documentclass[sn-mathphys,Numbered]{sn-jnl}

\usepackage{graphicx}%
\usepackage{multirow}%
\usepackage{amsmath,amssymb,amsfonts}%
\usepackage{amsthm}%
\usepackage{mathrsfs}%
\usepackage[title]{appendix}%
\usepackage{xcolor}%
\usepackage{textcomp}%
\usepackage{manyfoot}%
\usepackage{booktabs}%
\usepackage{algorithm}%
\usepackage{algorithmicx}%
\usepackage{algpseudocode}%
\usepackage{listings}%
\usepackage{enumerate}
\usepackage{setspace}

\theoremstyle{thmstyleone}%
\newtheorem{theorem}{Theorem}[section]%
\newtheorem{lemma}{Lemma}[section]%

\theoremstyle{thmstyletwo}%

\theoremstyle{thmstylethree}%
\newtheorem{definition}{Definition}[section]%

\theoremstyle{thmstylefour}%
\newtheorem{corollary}{Corollary}[section]%

\numberwithin{equation}{section}
\newcommand{\Om}{\Omega}
\newcommand{\cA}{\mathcal A}
\newcommand{\cL}{\mathcal L}

\newcommand{\divg}{\operatorname{div}_{g}}
\newcommand{\grad}{\nabla}
\newcommand{\dvol}{\,dV_g}

\newcommand{\norm}[1]{\left\lVert #1\right\rVert}
\newcommand{\abs}[1]{\left\lvert #1\right\rvert}
\newcommand{\ip}[2]{\left\langle #1,#2\right\rangle}

\begin{document}
	
	\title[Strong unique continuation for a Riemannian degenerate elliptic operator]{Strong unique continuation on Riemannian manifolds for degenerate elliptic operators with convexified weights}
	
	
	\author*[1]{\fnm{Bin} \sur{Shen}}\email{shenbin@seu.edu.cn}
	
	\author[2]{\fnm{Zisu} \sur{Zhao}}\email{steven1998zhao@outlook.com}
	
	\equalcont{These authors contributed equally to this work.}
	
	\affil*[1]{\orgdiv{School of Mathematics}, \orgname{Southeast University}, \orgaddress{\street{Dongnandaxue Road 2}, \city{Nanjing}, \postcode{211189}, \state{Jiangsu}, \country{China}}}
	\affil*[2]{\orgdiv{Institute of Mathematics}, \orgname{ Academy of Mathematics and Systems Science}, \orgname{Chinese Academy of Sciences}, \orgaddress{\street{Zhongguancun East Road 55}, \postcode{100190}, \state{Beijing}, \country{China}}}


	
	\abstract{
				In this manuscript, we consider a degenerate elliptic equation on a Riemannian manifold whose weight vanishes at an interior point $p$, and prove the strong unique continuation property (SUCP) at $p$. The curved setting gives rise to intrinsic geometric error terms that obstruct the analytic estimates available in the Euclidean case. To overcome this, we introduce a convexified radial weight that restores coercivity in the Carleman commutator. Combining the resulting Carleman estimate with a Hardy-type inequality and a Caccioppoli inequality, we show that any weak solution vanishing to infinite order at the degeneracy point must be identically zero on the connected domain. This result is the first SUCP theorem for this class of degenerate operators on general Riemannian manifolds.}

	\keywords{Strong unique continuation, Degenerate elliptic operator, Carleman estimate, Hardy inequality, Convexification, Riemannian manifold}
	
	
	\pacs[MSC Classification]{58J05, 35J70, 35B60}
	
	\maketitle

	\section{Introduction}
	
	The study of degenerate elliptic operators on Riemannian manifolds is motivated by a range of problems in geometric analysis, spectral theory, and mathematical physics. The operator
\begin{equation}\label{eq:main-equation}
    \mathcal{L}_{\alpha,V}u
    :=
    -\operatorname{div}_g(r^\alpha \nabla u)+Vu,
    \qquad 0<\alpha<2,
\end{equation}
where $r(x)=d_g(p,x)$ denotes the geodesic distance to a fixed interior point $p\in M$, arises naturally in several such problems.

First, \eqref{eq:main-equation} arises in the analysis of elliptic operators on spaces with conical or edge singularities. Let
\[
C(Y)=(0,\infty)_r\times Y
\]
be a metric cone over a closed Riemannian manifold $(Y^{n-1},h)$, equipped with the conic metric
\[
g=dr^2+r^2h.
\]
The scalar Laplacian on $C(Y)$ takes the form
\[
\Delta_{C(Y)}
=
-\partial_r^2-\frac{n-1}{r}\partial_r+\frac{1}{r^2}\Delta_Y,
\]
which, after multiplication by $r^\alpha$, leads to an operator of the type \eqref{eq:main-equation}, with $p$ corresponding to the cone tip. More generally, on manifolds with edge singularities, where the singular stratum has dimension $k\geq 1$, the model operator near the edge is of the form
\[
-\operatorname{div}_g(r^\alpha\nabla u)
+\text{lower-order terms},
\]
where $r$ denotes the distance to the edge and $\alpha$ encodes the asymptotic geometry of the cross-section \cite{Mazzeo1991}.

The systematic study of such operators was initiated by Mazzeo \cite{Mazzeo1991} and further developed by Mazzeo and Vertman \cite{MazzeoVertman2014} in the context of edge-degenerate operators. These works established Fredholm theory, heat-kernel asymptotics, and index theory for operators with polyhomogeneous coefficients on singular spaces. In particular, the linearized stability analysis of geometric flows, such as the Ricci flow or mean curvature flow, near singularities gives rise to operators of the form \eqref{eq:main-equation}, where the exponent $\alpha$ is determined by the asymptotic decay rate of the background metric \cite{BahuaudVertman2016}. Thus, understanding the qualitative properties of solutions to \eqref{eq:main-equation}, especially their behavior near the degeneracy, is an important prerequisite for nonlinear analysis on singular spaces.

The second, independent motivation comes from the theory of weighted Hardy inequalities on Riemannian manifolds. In Euclidean space, the classical Hardy inequality
\[
\int_{\mathbb{R}^n}\frac{|u|^2}{|x|^2}\,dx
\leq
\left(\frac{2}{n-2}\right)^2
\int_{\mathbb{R}^n}|\nabla u|^2\,dx,
\qquad n\geq 3,
\]
has a natural weighted generalization involving $|x|^\alpha$. The critical exponent $\alpha=2$ corresponds to the sharp constant $1/4$, whereas $0<\alpha<2$ places the operator in a case where the Hardy weight $r^{\alpha-2}$ is relatively compact with respect to the Dirichlet form $\int r^\alpha|\nabla u|^2$, see \cite{Davies1999}.

On a Riemannian manifold, Hardy-type inequalities naturally reflect the local geometry of geodesic spheres. Results such as those of Carron show that, under suitable geometric hypotheses, curvature and volume growth can influence the form and sharpness of Hardy inequalities on complete manifolds \cite{Carron2004}. In the present local setting, the relevant geometric quantity is the radial volume defect
\[
    h:=r\Delta_g r-(n-1),
\]
which measures the deviation of the radial volume density from its Euclidean counterpart \cite{Petersen2016}. It appears in the exact divergence identity
\[
    \operatorname{div}_g(X)=n+h,
    \qquad
    X=r\nabla r,
    \tag{1.2}
\]
and therefore enters both the weighted Hardy inequality and the subsequent Carleman calculation. In Euclidean space one has $h\equiv0$, whereas on a general smooth Riemannian manifold $h=O(r^2)$ in a sufficiently small normal neighborhood of $p$. Thus curvature appears in the argument through local geometric error terms, rather than through an imposed curvature sign condition. Since the analysis is carried out in a sufficiently small normal ball, these errors can be controlled by localization, and no additional curvature assumption is required for the qualitative unique continuation result. In this sense, the operator \eqref{eq:main-equation} provides a natural setting in which the effects of geometric variation and pointwise degeneracy can be studied simultaneously.

Moreover, \eqref{eq:main-equation} is closely related to the classical theory of the strong unique continuation property (SUCP) and its extension to degenerate operators. The SUCP asserts that any solution vanishing to infinite order at a point must vanish identically. For uniformly elliptic operators with smooth coefficients, this property is well established. Garofalo and Lin introduced the frequency-function method, which relates the vanishing order of a solution to the growth of its Dirichlet energy \cite{GarofaloLin1987}. Independently, Koch and Tataru developed robust Carleman estimates for operators with nonsmooth coefficients, obtaining almost optimal regularity conditions for the SUCP \cite{KochTataru2001}. It is also worth mentioning Rüland's work on the Caffarelli--Silvestre extension over a compact Riemannian manifold \cite{Ruland2017}. Although that work also concerns degenerate elliptic operators, the degeneracy occurs at the boundary rather than at an interior point.

However, when the operator degenerates at an interior point, as in \eqref{eq:main-equation}, both approaches encounter fundamental obstacles. The frequency-function method relies on the monotonicity of the frequency function, which in turn depends on an $L^2$-energy identity involving the vector field $X=r\nabla r$.

On a manifold, the geometric error terms $h$ and $B$ (see \eqref{eq:h-B} below) disrupt the exact differential identities underlying the monotonicity formula. To the best of our knowledge, no satisfactory proof based on the frequency-function method is available for the interior degenerate case.

Carleman estimates, on the other hand, are more flexible, but also encounter the same difficulty as frequency method. In Euclidean space, the radial vector field $X=x\cdot\nabla$ satisfies the exact identities
\begin{equation}\label{2.3}
    \operatorname{div}(X)=n,
    \qquad
    \nabla X=g.
\end{equation}
These identities yield exact cancellations in the Carleman commutator computation. More precisely, the terms involving derivatives of the weight function either vanish or combine into a manifestly positive expression. The resulting Carleman estimate takes a clean form, with the two leading terms
\[
\tau^3\int \rho^{-2\lambda-4}u^2
\qquad\text{and}\qquad
\tau\int \rho^{-2\lambda-2}r^\alpha|\nabla u|^2
\]
controlled by the $L^2$-norm of the operator applied to $u$ \cite{YangWuGuo2026}.

On a general Riemannian manifold, however, the identities in \eqref{2.3} are replaced by
\begin{equation}
    \operatorname{div}_gX=n+h,
    \qquad
    \nabla X=g+B,
\end{equation}
where
\begin{equation}\label{eq:h-B}
    h:=r\Delta_g r-(n-1),
    \qquad
    B:=\nabla X-g.
\end{equation}
The geometric meaning of the error terms in \eqref{eq:h-B} is worth noting. The scalar $h$ measures the deviation of the volume growth of geodesic spheres from its Euclidean counterpart. More precisely, in geodesic polar coordinates, the Riemannian volume element near $p$ can be written as
\[
dV_g=r^{n-1}q(r,\theta)\,dr\,d\sigma,
\qquad
q(r,\theta)=1+O(r^2),
\]
and
\[
h=r\partial_r\log q.
\]
The tensor $B$ measures the deviation of the radial vector field from the Euclidean homothetic identity $\nabla X=g$. Equivalently,
\[
B=\operatorname{Hess}_g\left(\frac{r^2}{2}\right)-g.
\]
Thus, $B$ records the difference between the second fundamental form of geodesic spheres and its Euclidean model; in particular, it vanishes in Euclidean space \cite{Petersen2016}.

A standard expansion in geodesic normal coordinates, as stated in Lemma \ref{lem:geometry}, gives estimates for $h$ and $B$ on a sufficiently small geodesic ball $B_R(p)$. In particular, both $h$ and $B$ vanish to second order at $p$. One might therefore expect that these terms could be treated as small perturbations and the Euclidean argument recovered by taking $R$ sufficiently small. However, this observation alone is not sufficient. In the Carleman commutator, $h$ and $B$ do not occur merely as additive perturbations. They are accompanied by the large Carleman parameter $\tau$, which tends to infinity. Although $r^2$ is small on a sufficiently small ball, the product $\tau r^2$ is not uniformly small as $\tau\to\infty$. Hence, shrinking $R$ alone does not compensate for the growth of $\tau$.


To handle this obstruction, we introduce the convexified radial weight
\begin{equation}\label{eq:Phi}
    \Phi(r):=-\log\rho(r)+\rho(r)^\epsilon,
    \qquad
    0<\epsilon<1,
\end{equation}
where
\begin{equation}
    \rho:=r^\gamma,
    \qquad
    \gamma:=1-\frac{\alpha}{2}\in(0,1).
\end{equation}
This choice is dictated by the structure of the Carleman commutator. For the weight \eqref{eq:Phi}, the associated quantities
\begin{equation}\label{eq:a-b}
    a:=X\Phi=-\gamma\bigl(1-\epsilon\rho^\epsilon\bigr),
    \qquad
    b:=Xa=\gamma^2\epsilon^2\rho^\epsilon
\end{equation}
satisfy $b>0$. By contrast, for the pure logarithmic weight, corresponding to $\epsilon=0$, one has
\[
a=-\gamma,
\qquad
b=0.
\]
Consequently, the additional positive contribution generated by $b$ is absent in the purely logarithmic case, and the remaining positive terms are insufficient to control the geometric errors involving $h$ and $B$. The convexification term $\rho^\epsilon$ produces a positive contribution of the form
\begin{equation*}
    \tau^3\int \rho^\epsilon z^2\,dV_g,
\end{equation*}
which is strong enough to absorb the relevant terms involving $h$, $B$, and their radial derivatives. At the same time, since $\rho^\epsilon\to0$ as $r\to0$, the convexification is a lower-order modification of the logarithmic weight near $p$ and does not interfere with the Hardy or Caccioppoli estimates used in the final cutoff argument.

The convexification mechanism is therefore not merely a technical modification of the Euclidean proof. It is introduced to make corrections for the additional geometric terms that arise from the failure of the identities on manifolds in \eqref{2.3}. The purpose of this modification is to handle the geometric errors caused by curvature. In the Euclidean case, the identities in \eqref{2.3} give exact cancellations. On a Riemannian manifold, these cancellations leave small error terms involving \(h\) and \(B\). The extra term in the weight creates additional positive contributions in the Carleman estimate, which are used to control these errors. Under the local geometric bounds assumed below, this allows the commutator argument to be carried out without imposing additional curvature conditions. The same idea may also be useful in the study of other degenerate elliptic problems on Riemannian manifolds.

	Let $\Omega\subset M$ be a connected open set containing $p$, and define the weighted Sobolev space	
	\begin{equation}
		\mathcal{H}_\alpha^1(\Omega) := \left\{ u\in L^2(\Omega,dV_g) : u\in W_{\mathrm{loc}}^{1,2}(\Omega\setminus\{p\}),\ \int_\Omega r^\alpha|\nabla u|_g^2\,dV_g < \infty \right\},
	\end{equation}	
	with norm
	\begin{equation}
		\|u\|_{\mathcal{H}_\alpha^1(\Omega)}^2 := \int_\Omega u^2\,dV_g + \int_\Omega r^\alpha|\nabla u|_g^2\,dV_g.
	\end{equation}	
	A function $u\in \mathcal{H}_\alpha^1(\Omega)$ is a weak solution of \eqref{eq:main-equation} if
	\begin{equation}\label{eq:weak-formulation}
	    \int_\Omega r^\alpha \langle \nabla u,\nabla\phi\rangle_g\,dV_g + \int_\Omega Vu\phi\,dV_g = 0,	
	\end{equation}
	for every $\phi\in C_c^\infty(\Omega\setminus\{p\})$.
	\begin{definition}\label{def:infinite-order}
		A function $u\in L^2_{\mathrm{loc}}(\Omega,dV_g)$ is said to  {vanish to infinite order} at $p$ if, for every $k>0$,
		\begin{equation}\label{eq:infinite-order}
			\lim_{s\rightarrow 0} s^{-k} \int_{B_s(p)} u^2\,dV_g = 0.
		\end{equation}
	\end{definition}
	We prove the following result.
	\begin{theorem}[Main theorem]\label{thm:main}
		Let $u\in \mathcal{H}_\alpha^1(\Omega)$ be a weak solution of \eqref{eq:main-equation} with $V\in L^\infty(\Omega)$. If $u$ vanishes to infinite order at $p$ in the sense of \eqref{eq:infinite-order}, then $u\equiv 0$ on the entire connected domain $\Omega$.
	\end{theorem}
	
	This theorem provides the  {first SUCP result} for the operator \eqref{eq:main-equation} in the general Riemannian setting, extending the Euclidean result of Yang, Wu, and Guo \cite{YangWuGuo2026} to curved spaces without any restrictive curvature assumptions. 
	It demonstrates that the convexification of the weight is a useful method for absorbing geometric error terms in Carleman estimates, making it a technique that could potentially be applied to other problems involving degenerate operators on manifolds, such as the SUCP for operators with double degeneracies or for systems with critical Hardy potentials. Furthermore, the theorem confirms that the geometric structure alone (encoded in the boundedness of $h$ and $B$) is sufficient to guarantee SUCP for this class of operators, providing a positive result in a setting where counterexamples are known to exist for other types of degeneracies.\\

	The remainder of the paper is organized as follows.
	
	In Section 2, we collect the necessary geometric preliminaries. Lemma \ref{lem:geometry} records the estimates on $h$, $B$, and their radial derivatives that follow from the geodesic normal coordinate expansion. Lemma \ref{lem:hardy} proves a local Hardy-type inequality that accounts for the volume defect $h$, providing the functional-analytic foundation for the Carleman estimate. Lemma \ref{lem:caccioppoli} establishes the Caccioppoli inequality, which is used in the final cutoff argument.
	
	In Section 3, we present the Carleman estimate. After defining the convexified weight $\Phi$, we perform the commutator computation and prove the key estimate (Theorem \ref{thm:carleman}). The geometric error terms from $h$ and $B$ are systematically absorbed using the positive term $\tau^3\int \rho^\epsilon z^2\,dV_g$. A corollary with the potential $V$ is also stated.
	
	In Section 4, we complete the proof of Theorem \ref{thm:main}. Using a radial cutoff function $\phi_\delta$, we apply the Carleman estimate to $u_\delta = \phi_\delta u$, pass to the limit $\delta\to 0$ using the infinite-order vanishing assumption, and conclude that $u\equiv 0$ in a sufficiently small geodesic ball. The conclusion is then extended to the whole connected domain $\Omega$ via the standard weak unique continuation property away from the degeneracy.

	\section{Hardy-type and Caccioppoli-type inequalities}
	In this section, we collect the preliminary estimates that will be used throughout the proof of the Carleman estimate and the final cutoff argument. We begin with Lemma \ref{lem:geometry}, which records the growth estimates for the geometric error terms $h$ and $B$ (defined in \eqref{eq:h-B}) and their radial derivatives. These estimates follow from the standard geodesic normal coordinate expansion and are essential for controlling the terms that appear in the Carleman commutator. Based on these geometric bounds, we then establish a local Hardy-type inequality (Lemma \ref{lem:hardy}) that accounts for the Riemannian volume defect $h$. This inequality provides the weighted $L^2$-estimates needed to handle the critical exponent $q = n/\gamma$ as well as the subcritical exponents $\ell \neq q$. Finally, we recall the Caccioppoli inequality (Lemma \ref{lem:caccioppoli}), which gives a local energy estimate for weak solutions of \eqref{eq:main-equation} and will be used in the cutoff argument in Section 4 to localize the solution near the degeneracy.\\
	
	For the radial vector field $X = x\cdot \nabla$, we identify $\nabla X$ with the endomorphism $Y\mapsto\nabla_YX$; equivalently, $\nabla X(Y, Z)=g(\nabla_{Y}X, Z)$ for any given vector fields $Y$ and $Z$, so its associated symmetric $2$-tensor is $\operatorname{Hess}_g(r^2/2)$. By taking geodesic normal
	coordinates centered at $p$, the usual normal-coordinate expansion gives the
	following estimates.
	
	Choose $R_0>0$ such that
	\begin{equation}\label{eq:R0}
		0<R_0<\min\bigl\{1,\operatorname{inj}_g(p),
		\operatorname{dist}_g(p,M\setminus\Om)\bigr\}.
	\end{equation}
	All the estimates near $p$ will be carried out in $B_{R_0}(p)$.
	As defined in \eqref{eq:h-B}, the geometric terms $h$ and $B$ satisfy the following estimates.
	\begin{lemma}\label{lem:geometry}
		After reducing $R_0$ if necessary, there exists a constant $C>0$ such that
		on $B_{R_0}(p)\setminus\{p\}$,
		\begin{equation}\label{eq:geometry-bounds}
			\abs{h}+\abs{Xh}+\abs{X^2h}
			+r\abs{\nabla h}_g+\abs{B}_g+r\abs{\nabla B}_g
			\leq Cr^2.
		\end{equation}
		In particular,
		\begin{equation}\label{eq:div-X}
			\divg X=n+h.
		\end{equation}
	\end{lemma}
	
	\begin{proof}
		In geodesic normal coordinates, we have
		$g_{ij}=\delta_{ij}+O(r^2)$ and $\partial_kg_{ij}=O(r)$. The Gauss lemma gives
		$X=x^i\partial_i=r\nabla r$. Hence
		$\nabla X=g+O(r^2)$ and $\divg X=n+O(r^2)$. Differentiating these expansions
		in the radial direction gives \eqref{eq:geometry-bounds}.
	\end{proof}
	
	Next, we establish a Hardy-type inequality on Riemannian manifolds. The difference from
	the Euclidean case is the term $h$ in \eqref{eq:div-X}.
	
	\begin{lemma}[Hardy-type inequality]\label{lem:hardy}
		Let $0<R\leq R_0$ and
		$u\in C_0^{\infty}(B_R(p)\setminus\{p\})$. For $\ell\neq q$, set
		\[
		c_{\ell,R}:=\gamma\abs{q-\ell}
		-\norm{h}_{L^{\infty}(B_R(p))}.
		\]
		If $c_{\ell,R}>0$, then
		\begin{equation}\label{eq:hardy}
			\int_{B_R(p)}\rho^{-\ell}u^2\dvol
			\leq
			\left(\frac{2}{c_{\ell,R}}\right)^2
			\int_{B_R(p)}\rho^{-\ell}(Xu)^2\dvol.
		\end{equation}
		Moreover, put
		\[
		c_{q,R}:=\gamma-
		\norm{h\log\rho}_{L^{\infty}(B_R(p))}.
		\]
		If $c_{q,R}>0$, then
		\begin{equation}\label{eq:hardy-critical}
			\int_{B_R(p)}\rho^{-q}u^2\dvol
			\leq
			\left(\frac{2}{c_{q,R}}\right)^2
			\int_{B_R(p)}(\log\rho)^2\rho^{-q}(Xu)^2\dvol.
		\end{equation}
	\end{lemma}
	
	\begin{proof}
		It is clear that
		\[
		X\rho=\gamma\rho,
		\qquad
		X(\rho^{-\ell})=-\gamma\ell\rho^{-\ell}.
		\]
		It follows from \eqref{eq:div-X} that
		\begin{equation}\label{eq:hardy-div}
			\begin{aligned}
				\divg(\rho^{-\ell}X)
				&=X(\rho^{-\ell})+\rho^{-\ell}\divg X\\
				&=\bigl[\gamma(q-\ell)+h\bigr]\rho^{-\ell}.
			\end{aligned}
		\end{equation}
		Suppose first that $\ell\neq q$. Since
		$\abs{h}<\gamma\abs{q-\ell}$ on $B_R(p)$, the function
		$\gamma(q-\ell)+h$ has a fixed sign. Therefore,
		\begin{align*}
			c_{\ell,R}\int_{B_R(p)}\rho^{-\ell}u^2\dvol
			&\leq
			\left|\int_{B_R(p)}
			\divg(\rho^{-\ell}X)u^2\dvol\right|\\
			&=2\left|\int_{B_R(p)}\rho^{-\ell}uXu\dvol\right|\\
			&\leq2
			\left(\int_{B_R(p)}\rho^{-\ell}u^2\dvol\right)^{1/2}
			\left(\int_{B_R(p)}\rho^{-\ell}(Xu)^2\dvol\right)^{1/2}.
		\end{align*}
		This implies \eqref{eq:hardy}.
		
		For $\ell=q$, equation \eqref{eq:hardy-div} becomes
		\[
		\divg(\rho^{-q}X)=h\rho^{-q}.
		\]
		Since $X(\log\rho)=\gamma$, we have
		\begin{equation}\label{eq:critical-div}
			\divg\bigl(\rho^{-q}(\log\rho)X\bigr)
			=\bigl[\gamma+h\log\rho\bigr]\rho^{-q}.
		\end{equation}
		Using \eqref{eq:critical-div} and repeating the above computation, we obtain
		\begin{align*}
			c_{q,R}\int_{B_R(p)}\rho^{-q}u^2\dvol
			&\leq2\left|\int_{B_R(p)}
			\rho^{-q}(\log\rho)uXu\dvol\right|\\
			&\leq2
			\left(\int_{B_R(p)}\rho^{-q}u^2\dvol\right)^{1/2}
			\left(\int_{B_R(p)}(\log\rho)^2
			\rho^{-q}(Xu)^2\dvol\right)^{1/2}.
		\end{align*}
		Thus \eqref{eq:hardy-critical} holds.
	\end{proof}
	Such a Hardy-type inequality is also a generalization of the Euclidean result (cf. Lemma 2.1 in \cite{YangWuGuo2026}).
	By Lemma \ref{lem:geometry}, both
	$\norm{h}_{L^{\infty}(B_R(p))}$ and
	$\norm{h\log\rho}_{L^{\infty}(B_R(p))}$ tend to zero as $R\rightarrow0$.
	Hence the assumptions in Lemma \ref{lem:hardy} hold on a sufficiently small
	ball. \\
	

The Hardy inequality established above provides the functional-analytic control of the radial derivative $Xu$ in weighted $L^2$-spaces. This is precisely the estimate needed in the Carleman commutator computation in Section 3, where the radial derivative $X$ appears as the principal first-order operator. However, the Carleman estimate alone does not suffice to prove the strong unique continuation property. In the final cutoff argument, we also need a local energy estimate that controls the full gradient $\nabla u$ in terms of the $L^2$-norm of $u$ on a slightly larger ball. This is the role of the Caccioppoli inequality, which we record in the following subsection. Unlike the Hardy inequality, the Caccioppoli inequality does not require the delicate weight $\rho$ and follows from a standard cutoff argument using the equation itself.
	
	\begin{lemma}[Caccioppoli-type inequality]\label{lem:caccioppoli}
		There exists $C>0$, depending only on the local geometry and $\alpha$, such
		that every weak solution of \eqref{eq:main-equation} satisfies
		\begin{equation}\label{eq:caccioppoli}
			\int_{B_{(1-a)R}(p)}r^{\alpha}\abs{\nabla u}_g^2\dvol
			\leq
			\frac{C\bigl(1+\norm{V}_{L^{\infty}(\Om)}\bigr)}
			{a^2R^{2-\alpha}}
			\int_{B_R(p)}u^2\dvol,
		\end{equation}
		for $0<a<1$ and $0<R\leq R_0$.
	\end{lemma}
	
	\begin{proof}
		Choose a radial cutoff $\eta$ with $\eta=1$ on $B_{(1-a)R}(p)$,
		$\eta=0$ outside $B_R(p)$, and $\abs{\nabla\eta}\leq C/(aR)$.
		Using $\eta^2u$ as a test function in \eqref{eq:weak-formulation}, we obtain 
		\[
		\begin{split}
			0&=\int_{B_{R}(p)} <r^{\alpha}\nabla u,  \nabla \eta^{2}u>_{g} \dvol+ \int_{B_{R}(p)}  V \eta^{2}u^{2} \dvol\\
			&=\int_{B_{R}(p)} \eta^{2}|r^{\frac{\alpha}{2}}\nabla u|^{2}\dvol+2\int_{B_{R}(p)} \eta u r^{\alpha}<\nabla \eta, \nabla u>_{g}\dvol+\int_{B_{R}(p)} V \eta^{2}u^{2}\dvol.
		\end{split}
		\]		
		Since
		\begin{eqnarray*}
			\begin{split}
				-2\int_{B_{R}(p)} \eta u r^{\alpha}<\nabla \eta, \nabla u>_{g}\dvol\le &\frac{1}{2}\int_{B_{R}(p)} \left(\eta r^{\frac{\alpha}{2}}|\nabla u|\right)^{2}\dvol\\
				&+2\int_{B_{R}(p)} \left(r^{\frac{\alpha}{2}}|\nabla \eta| u\right)^{2}\dvol,
			\end{split}		
		\end{eqnarray*}	
		we have
		\begin{eqnarray*}
			\begin{split}
				\int_{B_{(1-a)R}(p)}r^{\alpha}\abs{\nabla u}_g^2\dvol&\leq \int_{B_{R}(p)} \eta^2r^{\alpha}\abs{\nabla u}^2\dvol \\
				&\leq
				4\int_{B_{R}(p)} r^{\alpha}\abs{\nabla\eta}^2u^2\dvol
				+2\norm{V}_{\infty}\int_{B_{R}(p)}\eta^2u^2\dvol.
			\end{split}
		\end{eqnarray*}
		On $B_R(p)$, $r^{\alpha}\leq R^{\alpha}$.  Since $R\leq1$ and
		$a\leq1$, both terms on the right are bounded by the right-hand side of
		\eqref{eq:caccioppoli}.
	\end{proof}

For later use, set
\[
    \cA_\alpha u:=\divg\bigl(r^\alpha\grad u\bigr).
\]
When applying the Carleman estimate to a cutoff of the solution, we use
the product formula
\[
    \cA_\alpha(\chi u)
    =\chi\cA_\alpha u
      +2r^\alpha\langle\grad\chi,\grad u\rangle_g
      +u\cA_\alpha\chi.
\]
Thus, in addition to the usual bound for \(\grad\chi\), we need an
estimate for \(\cA_\alpha\chi\).

Let \(0<s\leq R_0\), and consider the annulus
\[
    B_s(p)\setminus\overline{B_{s/2}(p)},
\]
where the distance \(r(x)=d_g(p,x)\) satisfies
\[
    \frac{s}{2}<r(x)<s.
\]
Suppose that \(\chi=\chi(r)\) is a radial cutoff whose derivatives satisfy
\[
    |\chi'(r)|\leq \frac{C}{s},
    \qquad
    |\chi''(r)|\leq \frac{C}{s^2}
\]
on this annulus. Here the primes denote differentiation with respect
to the radial variable.

Since \(|\grad r|_g=1\) in the punctured normal ball, the chain rule gives
\[
    \grad\chi=\chi'\grad r,
    \qquad
    \Delta_g\chi=\chi''+\chi'\Delta_g r.
\]
Consequently,
\[
\begin{aligned}
    \cA_\alpha\chi
    &=r^\alpha\Delta_g\chi
      +\langle\grad(r^\alpha),\grad\chi\rangle_g\\
    &=r^\alpha\chi''
      +\bigl(r^\alpha\Delta_g r+\alpha r^{\alpha-1}\bigr)\chi'\\
    &=r^\alpha\chi''
      +(n+\alpha-1+h)r^{\alpha-1}\chi',
\end{aligned}
\]
where the last equality follows from
\(r\Delta_g r=n-1+h\).

Lemma~\ref{lem:geometry} bounds \(h\) on the fixed normal ball.
Combining this bound with the estimates for \(\chi'\) and \(\chi''\),
we obtain
\[
    |\cA_\alpha\chi|
    \leq C\left(\frac{r^\alpha}{s^2}
                 +\frac{r^{\alpha-1}}{s}\right)
    =Cs^{\alpha-2}
      \left[
        \left(\frac rs\right)^\alpha
        +\left(\frac rs\right)^{\alpha-1}
      \right].
\]
Because \(1/2<r/s<1\), the expression in brackets is bounded by a
constant depending only on \(\alpha\). Hence
\begin{equation}\label{eq:cutoff-A}
    |\cA_\alpha\chi|\leq Cs^{\alpha-2}
    \qquad\text{on }B_s(p)\setminus\overline{B_{s/2}(p)}.
\end{equation}
The constant is independent of \(s\), provided that the constants in
the derivative bounds for \(\chi\) are fixed. In the final cutoff
argument, we apply this estimate separately on the inner and outer
annuli, taking \(s\) to be the corresponding radius.

	\section{Carleman estimate}
	In this section, we establish a Carleman estimate for the operator $\mathcal{A}_\alpha = \operatorname{div}_g(r^\alpha \nabla \cdot)$ using the convexified weight $\Phi = -\log \rho + \rho^\epsilon$ introduced in \eqref{eq:Phi}. The geometric error terms $h$ and $B$, controlled by Lemma \ref{lem:geometry}, are absorbed through the positive-definite commutator term generated by the $\rho^\epsilon$ correction. The main result (Theorem \ref{thm:carleman}) provides weighted $L^2$-estimates for functions supported away from $p$, with the radial derivative controlled by the Hardy inequality (Lemma \ref{lem:hardy}). A corollary incorporating the potential $V$ is also stated for later use in the cutoff argument.
	
We emphasize here that the Lemma \ref{lem:multiplier} is already the Carleman estimate, expressed for the
rescaled operator \(\mathscr L=\rho\mathcal A_\alpha\rho\). Theorem
\ref{thm:carleman} does not create a second Carleman mechanism, it translates the conjugated estimate back to \(u\) and recovers the exact gradient term used later.

Before proving Lemma~\ref{lem:multiplier}, we compare the Euclidean
calculation of Yang--Wu--Guo \cite{YangWuGuo2026} with its Riemannian
counterpart. In Euclidean space, the dilation field satisfies
\(\operatorname{div}X=n\) and \(\nabla X=g\), whereas on the manifold
\[
    \operatorname{div}_gX=n+h,
    \qquad
    \nabla X=g+B.
\]
The terms \(h\) and \(B\) therefore introduce corrections into the
integrations by parts used to estimate the cross terms. Although
\(|h|+|B|_g=O(r^2)\) near \(p\), these corrections have no prescribed
sign and occur with coefficients involving the large Carleman parameter.
The Euclidean cancellations alone therefore do not complete the argument.
We first identify these additional terms, which motivate the use of the
convexified weight in Lemma~\ref{lem:multiplier}.

Let \(0<R\leq R_0\), take
\(v\in C_0^\infty(B_R(p)\setminus\{p\})\), and put \(u=\rho^\beta v\).
A direct product calculation gives
\begin{equation}\label{eq:yang-product}
\begin{aligned}
    \mathcal A_\alpha u
    ={}&2\beta\gamma\rho^{\beta-2}Xv
      +\beta\gamma\bigl[\gamma(\beta+q-2)+h\bigr]
       \rho^{\beta-2}v
      +\rho^\beta\mathcal A_\alpha v.
\end{aligned}
\end{equation}
As in the Euclidean proof, choose \(\beta=(2\lambda+4-q)/2\). Then
\[
    -2\lambda+2\beta-4=-q,
    \qquad
    -2\lambda+2\beta-2=-q+2,
\]
which accounts for the weights in the two cross-term calculations below.
Since \(v\) and its derivatives vanish near both \(p\) and
\(\partial B_R(p)\), all integrations by parts have zero boundary
contributions.

Using \eqref{eq:hardy-div} with \(\ell=q\) and the divergence theorem,
we 

\begin{equation}\label{eq:first-error}
\begin{aligned}
    \int_{B_R(p)}\rho^{-q}X(v^2)\,dV_g
    &=-\int_{B_R(p)}
      \operatorname{div}_g(\rho^{-q}X)v^2\,dV_g\\
    &=-\int_{B_R(p)}h\rho^{-q}v^2\,dV_g.
\end{aligned}
\end{equation}
For the second cross term (see also Equation 2.14 in \cite{YangWuGuo2026}), integration by parts gives
\[
\begin{aligned}
    \int_{B_R(p)}\rho^{-q+2}(Xv)\mathcal A_\alpha v\,dV_g
    ={}&\gamma(q-2)\int_{B_R(p)}\rho^{-q}(Xv)^2\,dV_g\\
       &-\int_{B_R(p)}r^{2-n}
       \langle\nabla(Xv),\nabla v\rangle_g\,dV_g.
\end{aligned}
\]
Here we have used \(r^\alpha\rho^{-q+2}=r^{2-n}\).
Moreover, \(\nabla X=g+B\) implies
\[
    \langle\nabla(Xv),\nabla v\rangle_g
    =\frac12X(|\nabla v|_g^2)+|\nabla v|_g^2
      +B(\nabla v,\nabla v).
\]
Combining this identity with
\[
    \operatorname{div}_g(r^{2-n}X)=(2+h)r^{2-n}
\]
and integrating by parts once more, we find
\begin{equation}\label{eq:second-error}
\begin{aligned}
    \int_{B_R(p)}\rho^{-q+2}(Xv)\mathcal A_\alpha v\,dV_g
    ={}&\gamma(q-2)\int_{B_R(p)}\rho^{-q}(Xv)^2\,dV_g\\
       &+\int_{B_R(p)}r^{2-n}
       \left(\frac h2|\nabla v|_g^2
       -B(\nabla v,\nabla v)\right)\,dV_g.
\end{aligned}
\end{equation}

To control these geometric corrections, we use the weight
\begin{equation}\label{eq:Phi}
    \Phi=-\log\rho+\rho^\varepsilon,
    \qquad 0<\varepsilon<1.
\end{equation}
The additional term gives
\[
    X^2\Phi=\gamma^2\varepsilon^2\rho^\varepsilon>0,
\]
which supplies positive terms in the commutator calculation. The following
lemma shows that the geometric errors can be absorbed by these terms
when the radius is sufficiently small and the Carleman parameter is
sufficiently large.

\begin{lemma}[Carleman estimate for the conjugated operator]\label{lem:multiplier}
Fix \(0<\varepsilon<1\). There exist
\(R_\varepsilon\in(0,R_0]\), \(\tau_0>0\), and \(C>0\) such that, for every
\(z\in C_0^\infty(B_{R_\varepsilon}(p)\setminus\{p\})\) and
\(\tau\geq\tau_0\),
\begin{equation}\label{eq:multiplier}
\begin{aligned}
&\tau^3\int_{B_{R_\varepsilon}(p)}\rho^\varepsilon z^2\,dV_g
 +\tau\int_{B_{R_\varepsilon}(p)}
     \rho^\varepsilon r^2|\nabla z|_g^2\,dV_g\\
&\hspace{30mm}\leq
 C\int_{B_{R_\varepsilon}(p)}|\mathscr L_\tau z|^2\,dV_g,
\end{aligned}
\end{equation}
where
\[
    \mathscr L:=\rho\mathcal A_\alpha\rho,
    \qquad
    \mathscr L_\tau z:=e^{\tau\Phi}
       \mathscr L(e^{-\tau\Phi}z).
\]
\end{lemma}
	
\begin{proof}
We work on the common core
\[
    \mathscr D_{R_{\varepsilon}}
    :=C_0^\infty\bigl(B_{R_{\varepsilon}}(p)\setminus\{p\}\bigr)
    \subset L^2\bigl(B_{R_{\varepsilon}}(p),dV_g\bigr),
\]
so that every integration by parts below has no boundary contribution. All
formal adjoints are taken in this real Hilbert space. We write
\[
    (u,v)_{\varepsilon}
    :=\int_{B_{R_{\varepsilon}}(p)}uv\,dV_g,
    \qquad
    \|u\|_{\varepsilon}^2:=(u,u)_{\varepsilon}.
\]
Multiplication by a real-valued function is symmetric.

We first fix the tensor conventions. Let
\[
    \mathsf D_X(Y):=\nabla_YX,
    \qquad
    (\nabla X)(Y,Z):=g(\mathsf D_XY,Z).
\]
Thus $\mathsf D_X$ is an endomorphism, whereas $\nabla X$ is the associated
symmetric $(0,2)$-tensor. Lemma~\ref{lem:geometry} gives
\[
    \nabla X=g+B
\]
as an identity of $(0,2)$-tensors. We also write
\[
    \mathcal P:=I-\nabla r\otimes dr,
    \qquad
    \mathcal C:=r^2\mathcal P.
\]
Here $\mathcal P$ is the orthogonal projection onto the tangent spaces of the
geodesic spheres and $\mathcal C$ is a $(1,1)$-tensor. This notation avoids
confusing the tangential tensor with the generic constants used in the
estimates.

Put
\[
    \mathfrak m:=\operatorname{div}_gX=n+h.
\]
The geometric constant in Lemma~\ref{lem:geometry} will be denoted by
$C_{\mathrm{geo}}$.  We emphasize here that all constants below may
depend on $n,\alpha,\varepsilon$ and $C_{\mathrm{geo}}$.

We now turn to the operator calculation. The first task is to identify the
formal adjoints and separate the rescaled operator into its radial and
tangential parts. The divergence theorem applied to $uvX$ gives
\[
    \int_{B_{R_{\varepsilon}}(p)} (Xu)v\,dV_g
    =-\int_{B_{R_{\varepsilon}}(p)} u\bigl(Xv+\mathfrak m v\bigr)\,dV_g,
\]
and hence
\[
    X^*=-X-\mathfrak m.
\]
Define
\begin{equation}\label{eq:H-T}
    H:=X+\frac{\mathfrak m}{2},
    \qquad
    T:=\operatorname{div}_g(r^2\nabla\,\cdot)-X^2-\mathfrak m X.
\end{equation}
Then $H^*=-H$. Moreover,
\[
    \mathcal C\nabla z=r^2\nabla z-(Xz)X,
\]
so
\[
    Tz=\operatorname{div}_g(\mathcal C\nabla z).
\]
Thus $T^*=T$, and integration by parts gives
\begin{equation}\label{eq:T-energy}
\begin{aligned}
    -(Tz,z)_{\varepsilon}
    &=\int_{B_{R_{\varepsilon}}(p)}\langle\mathcal C\nabla z,\nabla z\rangle_g\,dV_g\\
    &=\int_{B_{R_{\varepsilon}}(p)}\bigl(r^2|\nabla z|_g^2-(Xz)^2\bigr)\,dV_g\geq0.
\end{aligned}
\end{equation}
The integrand is $r^2|\mathcal P\nabla z|_g^2$; hence $T$ measures exactly
the angular energy. If $f=f(r)$ is radial, then
$\mathcal C\nabla f=0$, and therefore
\[
    [T,f]=0.
\]

Let $\mathscr L:=\rho\cA_\alpha\rho$. Formula~\eqref{eq:yang-product} with
$\beta=1$ and the identity $\rho^2=r^{2-\alpha}$ give
\[
\begin{aligned}
    \mathscr Lz
    &=2\gamma Xz+\gamma\bigl[\gamma(q-1)+h\bigr]z
      +\rho^2\cA_\alpha z,\\
    \rho^2\cA_\alpha z
    &=\operatorname{div}_g(r^2\nabla z)-2\gamma Xz.
\end{aligned}
\]
The first-order terms cancel. Since
\[
    H^2
    =X^2+\mathfrak m X+\frac12Xh+\frac14(n+h)^2,
\]
we obtain
\begin{equation}\label{eq:L-form}
    \mathscr L=H^2+T+\kappa,
\end{equation}
where
\begin{equation}\label{eq:kappa}
    \kappa
    :=\gamma^2(q-1)+\gamma h
      -\frac12Xh-\frac14(n+h)^2.
\end{equation}
Thus $H^2$ is the radial symmetric part, $T$ is the angular symmetric part,
and $\kappa$ is a zeroth-order geometric coefficient. Lemma~\ref{lem:geometry}
now yields
\begin{equation}\label{eq:kappa-bound}
    |\kappa|\leq C_\kappa,
    \qquad
    |X\kappa|\leq C_\kappa r^2,
\end{equation}
with $C_\kappa=C_\kappa(n,\alpha,C_{\mathrm{geo}})$.

For the phase $\Phi=-\log\rho+\rho^\varepsilon$, set
\begin{equation}\label{eq:a-b}
    a:=X\Phi=-\gamma(1-\varepsilon\rho^\varepsilon),
    \qquad
    b:=Xa=\gamma^2\varepsilon^2\rho^\varepsilon.
\end{equation}
The positivity of $b=X^2\Phi$ is the strict convexity introduced by the
perturbation $\rho^\varepsilon$. After decreasing $R_\varepsilon$, we may
assume
\begin{equation}\label{eq:a-bound}
    \frac{\gamma}{2}\leq|a|\leq\gamma.
\end{equation}
Since $[H,f]=Xf$ for every multiplication operator $f$,
\[
    e^{\tau\Phi}He^{-\tau\Phi}=H-\tau a,
    \qquad
    (H-\tau a)^2=H^2-2\tau aH-\tau b+\tau^2a^2.
\]
The operator $T$ commutes with the radial weight. Consequently,
\begin{equation}\label{eq:conjugated}
\begin{aligned}
    \mathscr L_\tau z
    ={}&H^2z+Tz+\kappa z+\tau^2a^2z
       -2\tau aHz-\tau bz.
\end{aligned}
\end{equation}

We decompose the conjugated operator into its symmetric and
antisymmetric parts, following the standard approach to Carleman
estimates; see, for example, \cite[Section~2]{Yamamoto2009}, put
\[
\begin{aligned}
    S_\tau&:=H^2+T+\kappa+\tau^2a^2,\\
    A_\tau&:=-2\tau aH-\tau b.
\end{aligned}
\]
Then $S_\tau^*=S_\tau$. The term $-\tau b$ is precisely the correction that
makes $A_\tau$ antisymmetric: since $Ha=aH+b$,
\[
    A_\tau=-\tau(aH+Ha),
    \qquad
    A_\tau^*=-A_\tau.
\]
Thus $\mathscr L_\tau=S_\tau+A_\tau$, although $\mathscr L_\tau$ itself is
generally neither symmetric nor antisymmetric. Expanding its norm gives
\begin{equation}\label{eq:norm-identity}
\begin{aligned}
    \|\mathscr L_\tau z\|_{\varepsilon}^2
    ={}&\|S_\tau z\|_{\varepsilon}^2
       +\|A_\tau z\|_{\varepsilon}^2\\
       &+([S_\tau,A_\tau]z,z)_{\varepsilon}.
\end{aligned}
\end{equation}

Having written the conjugated operator as the sum of its symmetric and
antisymmetric parts, we turn to the commutator that provides the positive
terms.
Using $[H,f]=Xf$, $[T,a]=[T,b]=0$ and $H^*=-H$, one finds
\[
\begin{aligned}
    ([H^2,A_\tau]z,z)_{\varepsilon}
      &=4\tau\int_{B_{R_{\varepsilon}}(p)} b|Hz|^2\,dV_g
        -\tau\int_{B_{R_{\varepsilon}}(p)} X^2b\,z^2\,dV_g,\\
    [\tau^2a^2,A_\tau]&=4\tau^3a^2b,\\
    [T,A_\tau]&=-2\tau a[T,H],\\
    [\kappa,A_\tau]&=2\tau a(X\kappa).
\end{aligned}
\]
For example, the first identity follows from
\[
    [H^2,2aH+b]=4bH^2+4(Xb)H+X^2b
\]
and
\[
\begin{aligned}
    (bH^2z,z)_{\varepsilon}
    ={}&-\int_{B_{R_{\varepsilon}}(p)}b|Hz|^2\,dV_g-\int_{B_{R_{\varepsilon}}(p)}(Xb)zHz\,dV_g.
\end{aligned}
\]
The two mixed terms cancel. Adding the four contributions gives
\begin{equation}\label{eq:commutator}
\begin{aligned}
([S_\tau,A_\tau]z,z)_{\varepsilon}
={}&4\tau\int_{B_{R_{\varepsilon}}(p)} b|Hz|^2\,dV_g
   +4\tau^3\int_{B_{R_{\varepsilon}}(p)} a^2bz^2\,dV_g
   -\tau\int_{B_{R_{\varepsilon}}(p)} X^2b\,z^2\,dV_g\\
&-2\tau(a[T,H]z,z)_{\varepsilon}
   +2\tau\int_{B_{R_{\varepsilon}}(p)} a(X\kappa)z^2\,dV_g.
\end{aligned}
\end{equation}
The first two terms are the desired positive radial and zeroth-order
energies; the remaining terms are errors.

It remains to understand the geometric commutator. Since $T^*=T$ and
$H^*=-H$,
\[
    ([T,H]z,z)_{\varepsilon}
    =2(Tz,Hz)_{\varepsilon}.
\]
Writing $w=\nabla z$, integrating by parts, and using
\[
    \nabla(Xz)=\nabla_Xw+\mathsf D_X^*w,
    \qquad
    \nabla(\mathfrak m z)=\mathfrak m w+z\nabla\mathfrak m,
\]
together with
\[
    X\langle\mathcal Cw,w\rangle_g
    =\langle(\nabla_X\mathcal C)w,w\rangle_g
      +2\langle\mathcal C\nabla_Xw,w\rangle_g,
\]
we obtain
\begin{equation}\label{eq:TH-formula}
\begin{aligned}
([T,H]z,z)_{\varepsilon}
={}&\int_{B_{R_{\varepsilon}}(p)}\left\langle
   \bigl(\nabla_X\mathcal C
      -\mathcal C\circ\mathsf D_X
      -\mathsf D_X^*\circ\mathcal C\bigr)\nabla z,
      \nabla z\right\rangle_g\,dV_g\\
&-\int_{B_{R_{\varepsilon}}(p)} z\langle\mathcal C\nabla\mathfrak m,
      \nabla z\rangle_g\,dV_g.
\end{aligned}
\end{equation}
The $\mathfrak m$-terms cancel here because
\[
    \int_{B_{R_{\varepsilon}}(p)}Xf\,dV_g
    =-\int_{B_{R_{\varepsilon}}(p)}\mathfrak m f\,dV_g.
\]
Formula~\eqref{eq:TH-formula}
measures exactly the failure of radial and angular operators to commute; it
vanishes in the Euclidean model.

We now estimate it without changing the tensor identity $\nabla X=g+B$.
Because radial curves are geodesics,
\[
    \nabla_X\nabla r=0,
    \qquad
    \nabla_Xdr=0,
\]
and hence, after lowering one index of $\mathcal C$,
\[
    \nabla_X\mathcal C_{ij}
    =\nabla_X\bigl(r^2(g_{ij}-r_ir_j)\bigr)
    =2\mathcal C_{ij}.
\]
Using $\nabla_iX_j=g_{ij}+B_{ij}$, the two metric terms cancel:
\[
\begin{aligned}
\bigl(\nabla_X\mathcal C
   -\mathcal C\circ\mathsf D_X
   -\mathsf D_X^*\circ\mathcal C\bigr)_{ij}=-\mathcal C_{kj}g^{k\ell}B_{i\ell}
 -\mathcal C_{ik}g^{k\ell}B_{j\ell}.
\end{aligned}
\]
Since $\|\mathcal C\|_{\mathrm{op}}=r^2$, $|B|_g\leq C_{\mathrm{geo}}r^2$,
$\nabla\mathfrak m=\nabla h$, and $|\nabla h|_g\leq C_{\mathrm{geo}}r$, we get
\[
\begin{aligned}
&\left|\left\langle
   \bigl(\nabla_X\mathcal C
      -\mathcal C\mathsf D_X
      -\mathsf D_X^*\mathcal C\bigr)\nabla z,
   \nabla z\right\rangle_g\right|
   \leq C_{\mathrm{geo}}r^4|\nabla z|_g^2,\\
&|z\langle\mathcal C\nabla\mathfrak m,\nabla z\rangle_g|
   \leq C_{\mathrm{geo}}
      \bigl(r^4|\nabla z|_g^2+r^2z^2\bigr).
\end{aligned}
\]
The factor $a$ causes no new derivative term because it is radial and
$\mathcal C\nabla a=0$. Therefore
\begin{equation}\label{eq:TH-bound}
    |(a[T,H]z,z)_{\varepsilon}|
    \leq C_{TH}\int_{B_{R_{\varepsilon}}(p)}
       \bigl(r^4|\nabla z|_g^2+r^2z^2\bigr)\,dV_g,
\end{equation}
where $C_{TH}$ depends on $\gamma$ and $C_{\mathrm{geo}}$. The two additional
powers of $r$ are what make this error absorbable on a small ball.

The commutator controls $Hz$ and $z$, whereas the statement of the lemma
asks for the full gradient. We therefore compare $Hz$ with the radial
derivative $Xz$ and then recover the tangential part. Expanding
$H=X+\mathfrak m/2$ and integrating the mixed term with
$X^*=-X-\mathfrak m$ gives
\begin{equation}\label{eq:H-expansion}
\begin{aligned}
    \int_{B_{R_{\varepsilon}}(p)} b|Hz|^2\,dV_g
    ={}&\int_{B_{R_{\varepsilon}}(p)} b(Xz)^2\,dV_g\\
       &-\int_{B_{R_{\varepsilon}}(p)} b\left(
          \frac{\mathfrak m^2}{4}
          +\frac{\gamma\varepsilon\mathfrak m}{2}
          +\frac{X\mathfrak m}{2}
          \right)z^2\,dV_g.
\end{aligned}
\end{equation}
Indeed, $Xb=\gamma\varepsilon b$ and
\[
    \int_{B_{R_{\varepsilon}}(p)} b\mathfrak m zXz\,dV_g
    =-\frac12\int_{B_{R_{\varepsilon}}(p)} b
      \bigl(\gamma\varepsilon\mathfrak m
            +X\mathfrak m+\mathfrak m^2\bigr)z^2\,dV_g.
\]
Thus $Hz$ differs from the true radial derivative only by a zeroth-order
term. Applying Lemma~\ref{lem:hardy} with $\ell=-\varepsilon$ yields
\[
    \int_{B_{R_{\varepsilon}}(p)} bz^2\,dV_g
    \leq
    \frac{4}{\bigl(\gamma(q+\varepsilon)
       -\|h\|_{L^\infty(B_{R_{\varepsilon}}(p))}\bigr)^2}
    \int_{B_{R_{\varepsilon}}(p)} b(Xz)^2\,dV_g.
\]
As the provisional radius $R_{\varepsilon}\rightarrow0$, the coefficient
that must be absorbed tends to
\[
    \frac{n^2+2n\gamma\varepsilon}
         {(n+\gamma\varepsilon)^2}
    =1-\frac{\gamma^2\varepsilon^2}
       {(n+\gamma\varepsilon)^2}<1.
\]
The strict gap is Euclidean; curvature only determines how small the radius
must be through $C_{\mathrm{geo}}R_{\varepsilon}^2$. After decreasing
$R_\varepsilon$, we therefore obtain
\begin{equation}\label{eq:H-coercive}
    \int_{B_{R_{\varepsilon}}(p)} b|Hz|^2\,dV_g
    \geq c_\varepsilon\int_{B_{R_{\varepsilon}}(p)} b(Xz)^2\,dV_g.
\end{equation}

To recover the tangential derivatives, define
\begin{equation}\label{eq:G}
\begin{aligned}
    G[z]
    &:=-\bigl(Tz,bz\bigr)_{\varepsilon}\\
    &=\int_{B_{R_{\varepsilon}}(p)} b\bigl(r^2|\nabla z|_g^2-(Xz)^2\bigr)\,dV_g
      =\int_{B_{R_{\varepsilon}}(p)} br^2|\mathcal P\nabla z|_g^2\,dV_g\geq0.
\end{aligned}
\end{equation}
Here the derivative of $b$ disappears because $\mathcal C\nabla b=0$.
Using the definition of $S_\tau$ and $H^*=-H$, we have
\[
\begin{aligned}
    G[z]
    ={}&-(S_\tau z,bz)_{\varepsilon}
       -\int_{B_{R_{\varepsilon}}(p)} b|Hz|^2\,dV_g
       -\int_{B_{R_{\varepsilon}}(p)}(Xb)zHz\,dV_g\\
       &+\int_{B_{R_{\varepsilon}}(p)}(\kappa+\tau^2a^2)bz^2\,dV_g.
\end{aligned}
\]
After multiplication by $\tau$, the favorable negative term is discarded,
and Cauchy--Schwarz and Young's inequalities give, for every $\eta>0$,
\begin{equation}\label{eq:G-estimate}
\begin{aligned}
    \tau G[z]
    \leq{}&\eta\|S_\tau z\|_{\varepsilon}^2
       +\eta\tau\int_{B_{R_{\varepsilon}}(p)} b|Hz|^2\,dV_g\\
       &+C_\eta\tau^3\int_{B_{R_{\varepsilon}}(p)} a^2bz^2\,dV_g.
\end{aligned}
\end{equation}
We used $Xb=\gamma\varepsilon b$, $|\kappa|\leq C_\kappa$,
$0<b\leq C_0$, $|a|\geq\gamma/2$ and $\tau\geq1$. Notice that this estimate
comes directly from the definition of $S_\tau$; it does not use the desired
Carleman inequality.

Because $b=\gamma^2\varepsilon^2r^{\gamma\varepsilon}$, on
$B_{R_{\varepsilon}}(p)$ we have
\[
    r^2\leq d_{R_{\varepsilon}}b,
    \qquad
    d_{R_{\varepsilon}}:=\frac{R_{\varepsilon}^{2-\gamma\varepsilon}}
                 {\gamma^2\varepsilon^2}
    \longrightarrow0
    \quad\text{as }R_{\varepsilon}\rightarrow0.
\]
Combining this with \eqref{eq:TH-bound}, the radial--angular identity, and
\eqref{eq:G}, we obtain
\[
\begin{aligned}
    2\tau|(a[T,H]z,z)_{\varepsilon}|
    \leq C_{TH}d_{R_{\varepsilon}}\tau
    \left(
       \int_{B_{R_{\varepsilon}}(p)}b(Xz)^2\,dV_g
       +G[z]
       +\int_{B_{R_{\varepsilon}}(p)}bz^2\,dV_g
    \right).
\end{aligned}
\]
We use \eqref{eq:H-coercive} for the first term,
\eqref{eq:G-estimate} for the second, and
\[
    \tau\int_{B_{R_{\varepsilon}}(p)} bz^2\,dV_g
    \leq\frac{4}{\gamma^2}\tau^3
       \int_{B_{R_{\varepsilon}}(p)} a^2bz^2\,dV_g
\]
for the third. Fixing $\eta$ and then choosing $R_\varepsilon$ sufficiently
small gives
\begin{equation}\label{eq:TH-absorb}
\begin{aligned}
    2\tau|(a[T,H]z,z)_{\varepsilon}|
    \leq{}&\frac14\|S_\tau z\|_{\varepsilon}^2
       +\frac12\tau\int_{B_{R_{\varepsilon}}(p)} b|Hz|^2\,dV_g\\
       &+\frac12\tau^3\int_{B_{R_{\varepsilon}}(p)} a^2bz^2\,dV_g.
\end{aligned}
\end{equation}
This is the only absorption that requires the spatial ball to be small.

The remaining errors are zeroth order. Since
$X^2b=(\gamma\varepsilon)^2b$ and $|X\kappa|\leq C_\kappa r^2$, after
$R_\varepsilon$ has been fixed they are of order $\tau$, whereas the good
zeroth-order term in \eqref{eq:commutator} is of order $\tau^3$. Increasing
$\tau_0$ therefore gives
\begin{equation}\label{eq:lower-absorb}
    \tau\left|\int_{B_{R_{\varepsilon}}(p)} X^2b\,z^2\,dV_g\right|
    +2\tau\left|\int_{B_{R_{\varepsilon}}(p)} a(X\kappa)z^2\,dV_g\right|
    \leq\frac12\tau^3\int_{B_{R_{\varepsilon}}(p)} a^2bz^2\,dV_g.
\end{equation}

At this point all error terms have been absorbed. It remains only to
assemble the estimates and recover the exact norm appearing in the lemma.
Insert \eqref{eq:commutator}, \eqref{eq:TH-absorb} and
\eqref{eq:lower-absorb} into \eqref{eq:norm-identity}, and discard the
nonnegative term $\|A_\tau z\|_{\varepsilon}^2$. Using
$a^2\geq\gamma^2/4$ and weakening fixed numerical constants, we obtain
\begin{equation}\label{eq:basic-coercive}
\begin{aligned}
    \|\mathscr L_\tau z\|_{\varepsilon}^2
    \geq{}&c\tau\int_{B_{R_{\varepsilon}}(p)} b|Hz|^2\,dV_g
       +c\tau^3\int_{B_{R_{\varepsilon}}(p)} bz^2\,dV_g\\
       &+c\|S_\tau z\|_{\varepsilon}^2.
\end{aligned}
\end{equation}
This is the coercive core of the proof. It controls $Hz$, but the statement
of the lemma requires the full gradient. By \eqref{eq:H-coercive} and
\eqref{eq:basic-coercive},
\[
    \tau\int_{B_{R_{\varepsilon}}(p)} b(Xz)^2\,dV_g
    \leq C_0\|\mathscr L_\tau z\|_{\varepsilon}^2.
\]
Applying \eqref{eq:G-estimate} after \eqref{eq:basic-coercive} also gives
\[
    \tau G[z]\leq C_0\|\mathscr L_\tau z\|_{\varepsilon}^2.
\]
Since
\[
    \int_{B_{R_{\varepsilon}}(p)} br^2|\nabla z|_g^2\,dV_g
    =\int_{B_{R_{\varepsilon}}(p)} b(Xz)^2\,dV_g+G[z],
\]
we conclude that
\[
    \tau\int_{B_{R_{\varepsilon}}(p)} br^2|\nabla z|_g^2\,dV_g
    \leq C_0\|\mathscr L_\tau z\|_{\varepsilon}^2.
\]
Finally, $b=\gamma^2\varepsilon^2\rho^\varepsilon$, so this estimate and
\eqref{eq:basic-coercive} give exactly \eqref{eq:multiplier}. The order of
choices is $ \eta$  fixed,$R_{\varepsilon}$ chosen small, $\tau_{0}$ chosen large.
\end{proof}

	We can now prove the Carleman estimate. 
	
	\begin{theorem}\label{thm:carleman}
		Fix $0<\varepsilon<1$. There exist $R_{\varepsilon}\in(0,R_0]$,
		$C>0$, and $\lambda_0>0$ such that, for every
		$u\in C_0^{\infty}(B_{R_{\varepsilon}}(p)\setminus\{p\})$ and
		$\lambda\geq\lambda_0$,
		\begin{equation}\label{eq:carleman}
			\begin{aligned}
				&\lambda^2\int
				\rho^{-2\lambda-4+\varepsilon}
				e^{2(\lambda+1)\rho^{\varepsilon}}u^2\dvol+\int\rho^{-2\lambda-2+\varepsilon}e^{2(\lambda+1)\rho^{\varepsilon}}r^{\alpha}\abs{\nabla u}_g^2\dvol\\
				\leq&
				C\int
				\rho^{-2\lambda}e^{2(\lambda+1)\rho^{\varepsilon}}
				\abs{\cA_{\alpha}u}^2\dvol.
			\end{aligned}
		\end{equation}
	\end{theorem}
	
	\begin{proof}
		Put $\tau=\lambda+1$.
		
		\noindent\textbf{Step 1.} Let
		\[
		z=e^{\tau\Phi}\rho^{-1}u.
		\]
		Then
		\[
		\mathscr L_{\tau}z
		=e^{\tau\Phi}\rho\cA_{\alpha}u.
		\]
		Applying Lemma \ref{lem:multiplier}, and using
		$e^{2\tau\Phi}=\rho^{-2\tau}e^{2\tau\rho^{\varepsilon}}$, we obtain
		\begin{equation}\label{eq:step1-zero}
			\begin{aligned}
				\lambda^2\int
				\rho^{-2\lambda-4+\varepsilon}
				e^{2(\lambda+1)\rho^{\varepsilon}}u^2\dvol
				\leq C\int
				\rho^{-2\lambda}e^{2(\lambda+1)\rho^{\varepsilon}}
				\abs{\cA_{\alpha}u}^2\dvol.
			\end{aligned}
		\end{equation}
		Here we have discarded one positive power of $\tau$.
		
		\noindent\textbf{Step 2.} It is evident that
		\begin{equation}\label{eq:energy-identity}
			r^{\alpha}\abs{\nabla u}_g^2
			=\frac12\cA_{\alpha}(u^2)-u\cA_{\alpha}u.
		\end{equation}
		Set
		\[
		W:=\rho^{-2\lambda-2+\varepsilon}
		e^{2(\lambda+1)\rho^{\varepsilon}}.
		\]
		It follows from \eqref{eq:energy-identity} that
		\begin{equation}\label{eq:step2-start}
			\begin{aligned}
				\int Wr^{\alpha}\abs{\nabla u}_g^2\dvol
				=\frac12\int W\cA_{\alpha}(u^2)\dvol
				-\int Wu\cA_{\alpha}u\dvol.
			\end{aligned}
		\end{equation}
		Since $W$ is radial and
		$\abs{X\log W}\leq C\lambda$ on $B_{R_{\varepsilon}}(p)$, integration by
		parts and Young's inequality give
		\begin{equation}\label{eq:step2-first}
			\begin{aligned}
				\frac12\left|\int W\cA_{\alpha}(u^2)\dvol\right|
				\leq{}&\frac14\int Wr^{\alpha}\abs{\nabla u}_g^2\dvol+C\lambda^2\int W\rho^{-2}u^2\dvol.
			\end{aligned}
		\end{equation}
		On the other hand,
		\begin{equation}\label{eq:step2-second}
			\begin{aligned}
				\left|\int Wu\cA_{\alpha}u\dvol\right|
				\leq{}&\frac14\lambda^2\int W\rho^{-2}u^2\dvol+\frac{C}{\lambda^2}\int W\rho^2
				\abs{\cA_{\alpha}u}^2\dvol.
			\end{aligned}
		\end{equation}
		Since $\rho^{\varepsilon}\leq1$ on a sufficiently small ball,
		$W\rho^2\leq\rho^{-2\lambda}e^{2(\lambda+1)\rho^{\varepsilon}}$.
		Combining \eqref{eq:step2-start}--\eqref{eq:step2-second}, and then applying
		\eqref{eq:step1-zero}, proves \eqref{eq:carleman}.
	\end{proof}
	
	Since $\cA_{\alpha}u=Vu-\cL_{\alpha,V}u$, Theorem
	\ref{thm:carleman} gives an additional term
	\[
	C\norm{V}_{\infty}^2
	\int\rho^{-2\lambda}e^{2(\lambda+1)\rho^{\varepsilon}}u^2\dvol.
	\]
	On a sufficiently small ball,
	$\rho^{-2\lambda}\leq\rho^{-2\lambda-4+\varepsilon}$.
	Thus this term can be absorbed into the first term on the left by choosing
	$\lambda$ sufficiently large. Therefore, we immediately have the following result.
	
	\begin{corollary}\label{cor:potential}
		There exist $C>0$ and $\lambda_1>0$ such that
		\begin{equation}\label{eq:carleman-potential}
			\begin{aligned}
				&\lambda^2\int
				\rho^{-2\lambda-4+\varepsilon}
				e^{2(\lambda+1)\rho^{\varepsilon}}u^2\dvol+
	                         \int
				\rho^{-2\lambda-2+\varepsilon}
				e^{2(\lambda+1)\rho^{\varepsilon}}
				r^{\alpha}\abs{\nabla u}_g^2\dvol\\
				&\hspace{25mm}\leq
				C\int
				\rho^{-2\lambda}e^{2(\lambda+1)\rho^{\varepsilon}}
				\abs{\cL_{\alpha,V}u}^2\dvol
			\end{aligned}
		\end{equation}
		holds for every compactly supported smooth function in the punctured ball,
		provided
		$$
		\lambda\geq\lambda_1
		\bigl(1+\norm{V}_{L^{\infty}(\Om)}\bigr).
		$$
	\end{corollary}
	
	
	\section{Strong unique continuation property}

We now combine the ingredients developed in the preceding sections to prove the main theorem. The argument begins by applying the Carleman estimate in Corollary~\ref{cor:potential} to the cutoff functions
\[
    u_\delta=\phi_\delta u,
\]
where $\phi_\delta$ removes a small neighborhood of the degenerate point $p$. The derivatives of the cutoff are supported in two annular regions: an inner annulus
$
    \{\delta\le r\le 2\delta\},
$
which shrinks to $p$ as $\delta\to0$, and a fixed outer annulus, for instance
$
    \left\{\frac{R}{2}\le r\le R\right\}.
$
The infinite-order vanishing assumption \eqref{eq:infinite-order} shows that the contribution from the inner annulus tends to zero as $\delta\to0$. After passing to the limit, we obtain an estimate on the fixed ball $B_{R/2}(p)$ in which the weighted $L^2$-norm of $u$ is controlled by terms supported in the outer annulus. Letting the Carleman parameter $\lambda\to\infty$, the separation of the weight between $B_{R/4}(p)$ and the outer annulus forces $u$ to vanish identically on $B_{R/4}(p)$. This proves strong unique continuation at the degenerate point. Finally, since the equation is uniformly elliptic with locally Lipschitz coefficients away from $p$, the classical weak unique continuation property propagates this local vanishing throughout the connected domain $\Omega$.

	\begin{proof}[Proof of Theorem \ref{thm:main}]
		Fix $0<\varepsilon<1$ and choose $0<R<R_{\varepsilon}$ so small that
		$\varepsilon\rho(R)^{\varepsilon}\leq1/2$. Then the function
		$\Phi=-\log\rho+\rho^{\varepsilon}$ is strictly decreasing on $(0,R)$.
		
		For $0<\delta<R/4$, choose a radial cutoff $\varphi_{\delta}$ satisfying
		\[
		\begin{array}{lll}
			\varphi_{\delta}=0&\text{when }r\leq\delta/2\text{ or }r\geq R,
			&\varphi_{\delta}=1\text{ when }\delta\leq r\leq R/2,
		\end{array}
		\]
		and
		\[
		\begin{array}{ll}
			\abs{\nabla\varphi_{\delta}}\leq C\delta^{-1},\quad
			\abs{\cA_{\alpha}\varphi_{\delta}}\leq C\delta^{\alpha-2},
			&\delta/2<r<\delta,\\[1mm]
			\abs{\nabla\varphi_{\delta}}\leq CR^{-1},\quad
			\abs{\cA_{\alpha}\varphi_{\delta}}\leq CR^{\alpha-2},
			&R/2<r<R.
		\end{array}
		\]
		Away from $p$, the equation is uniformly elliptic. By the standard interior
		regularity and a density argument, we may apply Corollary \ref{cor:potential}
		to $u_{\delta}:=\varphi_{\delta}u$. Since
		$\cL_{\alpha,V}u=0$, we have
		\begin{equation}\label{eq:commutator-cutoff}
			\cL_{\alpha,V}u_{\delta}
			=-2r^{\alpha}\ip{\nabla u}{\nabla\varphi_{\delta}}_g
			-u\cA_{\alpha}\varphi_{\delta}.
		\end{equation}
		Therefore,
		\begin{equation}\label{eq:cutoff-carl}
			\begin{aligned}
				&\lambda^2\int
				\rho^{-2\lambda-4+\varepsilon}
				e^{2(\lambda+1)\rho^{\varepsilon}}
				\varphi_{\delta}^2u^2\dvol\\
				&\quad\leq C\int
				\rho^{-2\lambda}e^{2(\lambda+1)\rho^{\varepsilon}}
				\left(r^{2\alpha}\abs{\nabla u}_g^2
				\abs{\nabla\varphi_{\delta}}_g^2
				+u^2\abs{\cA_{\alpha}\varphi_{\delta}}^2\right)\dvol.
			\end{aligned}
		\end{equation}
		
		We first deal with the inner annulus. For a fixed $\lambda$, using the cutoff
		bounds, Lemma \ref{lem:caccioppoli}, and \eqref{eq:infinite-order}, we obtain
		\begin{align*}
			&\int_{B_{\delta}(p)\setminus B_{\delta/2}(p)}
			\rho^{-2\lambda}e^{2(\lambda+1)\rho^{\varepsilon}}
			\left(r^{2\alpha}\abs{\nabla u}_g^2
			\abs{\nabla\varphi_{\delta}}_g^2
			+u^2\abs{\cA_{\alpha}\varphi_{\delta}}^2\right)\dvol\\
			&\qquad\leq
			C_{\lambda}\delta^{-2\lambda\gamma+2\alpha-4}
			\int_{B_{2\delta}(p)}u^2\dvol\longrightarrow0
		\end{align*}
		as $\delta\rightarrow0$. Thus, passing to the limit in
		\eqref{eq:cutoff-carl}, only the outer annulus remains. We have
		\begin{equation}\label{eq:outer-annulus}
			\begin{aligned}
				&\lambda^2\int_{B_{R/2}(p)}
				\rho^{-2+\varepsilon}e^{2(\lambda+1)\Phi}u^2\dvol\\
				&\quad\leq C_R\int_{B_R(p)\setminus B_{R/2}(p)}
				e^{2(\lambda+1)\Phi}
				\left(r^{\alpha}\abs{\nabla u}_g^2+u^2\right)\dvol.
			\end{aligned}
		\end{equation}
		Since $\Phi$ is decreasing,
		\[
		\Phi(r)\geq\Phi(R/4)\quad(r\leq R/4),
		\qquad
		\Phi(r)\leq\Phi(R/2)\quad(R/2\leq r\leq R).
		\]
		It follows from \eqref{eq:outer-annulus} that
		\begin{align*}
			\int_{B_{R/4}(p)}u^2\dvol
			\leq{}&\frac{C_R}{\lambda^2}
			e^{-2(\lambda+1)[\Phi(R/4)-\Phi(R/2)]}\\
			&\times\int_{B_R(p)\setminus B_{R/2}(p)}
			\left(r^{\alpha}\abs{\nabla u}_g^2+u^2\right)\dvol.
		\end{align*}
		Letting $\lambda\to\infty$, we conclude that
		$u\equiv0$ in $B_{R/4}(p)$.
		
		It remains to extend the conclusion to $\Om$. In a coordinate neighborhood
		separated from $p$, equation \eqref{eq:main-equation} can be written as
		\[
		-\partial_i(A^{ij}\partial_j u)+c(x)u=0,
		\qquad
		A^{ij}=\sqrt{\det g}\,r^{\alpha}g^{ij},
		\qquad
		c=\sqrt{\det g}\,V.
		\]
		The coefficients $A^{ij}$ are locally Lipschitz and uniformly elliptic on
		$\Om\setminus\{p\}$, and $c\in L^{\infty}$. The weak unique continuation
		property therefore implies $u\equiv0$ in the connected domain $\Om$.
	\end{proof}

	
	
	
	
	
	
	
%

\section*{Declarations}
\begin{itemize}
	\item Conflict of interest/Competing interests: All authors disclose no relevant relationships.
	\item Ethics approval and consent to participate: Not applicable.
	\item Data availability: Not applicable.
	\item Materials availability: Not applicable.
	\item Code availability: Not applicable.
	\item Generative AI and AI-assisted Technologies: This article does not use any generative AI assistance.
\end{itemize}
\bibliography{refs.bib}

\end{document}